%
\def\date{17 March 2014}  
\magnification=1200
\overfullrule=0pt
\newif\ifproofmode
\def\xrefsfilename{girth.xrf}  
\def\myinput#1{\immediate\openin0=#1\relax
   \ifeof0\write16{Cannot input file #1.}
   \else\closein0\input#1\fi}
\newcount\referno
\newcount\thmno
\newcount\secno
\referno=0
\thmno=0
\secno=0
\def\ifundefined#1{\expandafter\ifx\csname#1\endcsname\relax}
\myinput \xrefsfilename
\immediate\openout1=\xrefsfilename
\def\bibitem#1#2\par{\ifundefined{REFLABEL#1}\relax\else
 \global\advance\referno by 1\relax
 \immediate\write1{\noexpand\expandafter\noexpand\def
 \noexpand\csname REFLABEL#1\endcsname{\the\referno}}
 \global\expandafter\edef\csname REFLABEL#1\endcsname{\the\referno}
 \item{\the\referno.}#2\ifproofmode [#1]\fi\fi}
\def\cite#1{\ifundefined{REFLABEL#1}\ignorespaces
   \global\expandafter\edef\csname REFLABEL#1\endcsname{?}\ignorespaces
   \write16{ ***Undefined reference #1*** }\fi
 \csname REFLABEL#1\endcsname}
\def\nocite#1{\ifundefined{REFLABEL#1}\ignorespaces
   \global\expandafter\edef\csname REFLABEL#1\endcsname{?}\ignorespaces
   \write16{ ***Undefined reference #1*** }\fi}
\def\newthm#1#2\par{\global\advance\thmno by 1\relax
 \immediate\write1{\noexpand\expandafter\noexpand\def
 \noexpand\csname THMLABEL#1\endcsname{(\the\secno.\the\thmno)}}
 \global\expandafter\edef\csname THMLABEL#1\endcsname{(\the\secno.\the\thmno)}
 \bigbreak\penalty-800\noindent{\bf(\the\secno.\the\thmno)\enspace}\ignorespaces
 \ifproofmode {\bf[#1]} \fi{\sl#2}
 \medbreak\penalty-200}
\def\newsection#1#2\par{\global\advance\secno by 1\relax
 \immediate\write1{\noexpand\expandafter\noexpand\def
 \noexpand\csname SECLABEL#1\endcsname{\the\secno}}
 \global\expandafter\edef\csname SECLABEL#1\endcsname{\the\secno}
 \vskip0pt plus.3\vsize
 \vskip0pt plus-.3\vsize\bigskip\bigskip\vskip\parskip\penalty-250
 \message{\the\secno. #2}\thmno=0
 \centerline{\bf\the\secno. #2\ifproofmode {\rm[#1]} \fi}
 \nobreak\smallskip\noindent}
\def\refthm#1{\ifundefined{THMLABEL#1}\ignorespaces
 \global\expandafter\edef\csname THMLABEL#1\endcsname{(?)}\ignorespaces
 \write16{ ***Undefined theorem label #1*** }\fi
 \csname THMLABEL#1\endcsname}
\def\refsec#1{\ifundefined{SECLABEL#1}\ignorespaces
 \global\expandafter\edef\csname SECLABEL#1\endcsname{(?)}\ignorespaces
 \write16{ ***Undefined section label #1*** }\fi
 \csname SECLABEL#1\endcsname}

\font\smallrm=cmr8
\outer\def\thm#1#2\par{\medbreak\noindent{\bf(#1)\enspace}\ignorespaces
{\sl#2}\ifdim\lastskip<\medskipamount\removelastskip\penalty55\medskip\fi}
\def\cond#1#2\par{\smallbreak\noindent\rlap{\rm(#1)}\ignorespaces
\hangindent=36pt\hskip36pt{\rm#2}\smallskip}
\def\claim#1#2\par{{\medbreak\noindent\rlap{\rm(#1)}\ignorespaces
\rightskip20pt
\hangindent=20pt\hskip20pt{\ignorespaces\sl#2}\smallskip}}
\def\dfn#1{{\sl #1}}

\def\proof{\smallbreak\noindent{\sl Proof. }}



\def\qed{\hfill$\square$\bigskip\medskip}
\def\sqr#1#2{{\vcenter{\vbox{\hrule height.#2pt
\hbox{\vrule width.#2pt height #1pt \kern#1pt
\vrule width.#2pt}
\hrule height.#2pt}}}}
\def\square{\mathchoice\sqr56\sqr56\sqr{2.1}3\sqr{1.5}3}

\def\ref#1#2{\item{#1.}#2}
\outer\def\beginsection#1\par{\vskip0pt plus.3\vsize
   \vskip0pt plus-.3\vsize\bigskip\bigskip\vskip\parskip
   \message{#1}\centerline{\bf#1}\nobreak\smallskip\noindent}

\nopagenumbers
\footline={\hfil}
\baselineskip=12pt
\phantom{a}\vskip .25in
\centerline{\bf GIRTH SIX CUBIC GRAPHS}
\centerline{\bf HAVE PETERSEN MINORS}
\bigskip\bigskip
\centerline{Neil Robertson$^{*1}$\vfootnote{$^*$}{\smallrm Research partially
performed under a consulting agreement with Bellcore, and partially
supported by DIMACS Center,
Rutgers University, New Brunswick, New Jersey  08903, USA.}
\vfootnote{$^1$}{\smallrm Partially supported
by NSF under Grant No. DMS-9401981 and by ONR under Grant No.
N00014-92-J-1965.
}}
\centerline{Department of Mathematics}
\centerline{Ohio State University}
\centerline{231 W. 18th Ave.}
\centerline{Columbus, Ohio  43210, USA}
\bigskip
 
\centerline{P. D. Seymour}
\centerline{Bellcore}
\centerline{445 South St.}
\centerline{Morristown, New Jersey  07960, USA}
\bigskip
 
\centerline{and}
\bigskip
 
\centerline{Robin Thomas$^{*2}$\vfootnote{$^2$}{\smallrm
Partially supported
by NSF under Grants No. DMS-9623031 and DMS-1202640
and by ONR under Contract No.~N00014-93-1-0325.
}}
\centerline{School of Mathematics}
\centerline{Georgia Institute of Technology}
\centerline{Atlanta, Georgia  30332, USA}
\bigskip\bigskip\bigskip
\beginsection ABSTRACT

\parshape=1.5truein 5.5truein
We prove that every $3$-regular graph with no circuit of length less than
six has a subgraph isomorphic to a subdivision of the Petersen graph.

\vfill
\baselineskip 11pt
\noindent 3 January 1997\hfil\break
\noindent Revised \date
\vfil\eject
\baselineskip 18pt




\footline{\hss\tenrm\folio\hss}

\newsection{intro}INTRODUCTION

All \dfn{graphs} in this paper are finite, and may have loops and parallel
edges. A graph is \dfn{cubic} if the degree of every vertex (counting loops
twice) is three. The \dfn{girth} of a graph is the length of its shortest
circuit, or infinity if the graph has no circuits. (\dfn{Paths} and
\dfn{circuits} have no ``repeated" vertices.) The \dfn{Petersen graph} 
is the unique cubic graph of girth five on ten vertices. 
The Petersen graph is an obstruction
to many properties in graph theory, and often is, or is conjectured to be,
the {\sl only} obstruction. Such is the case for instance in the following
result of Alspach, Goddyn and Zhang [\cite{AlsGodZha}]. Let $G$ be
a graph, and let $p:E(G)\to{\bf Z}$ be a mapping. We say that $p$ is
\dfn{admissible} if $p(e)\ge0$ for every edge $e$ of $G$, and for every
edge-cut $C$, $\sum_{e\in C}p(e)$ is even and at least twice $p(f)$
for every edge $f\in C$.
We say that a graph $G$ is a \dfn{subdivision} of a graph $H$ if $G$
can be obtained from $H$ by replacing the edges of $H$ by internally
disjoint paths with the same ends and at least one edge.
We say that a graph
$G$ \dfn{contains} a graph $H$ if $G$ has a subgraph isomorphic to 
a subdivision of $H$.

\newthm{azg}For a graph $G$, the following two conditions are equivalent.
\item{{\rm(i)}}For every admissible mapping $p:E(G)\to{\bf Z}$ there 
exists a list of circuits of $G$ such that  every edge $e$ of $G$
belongs to precisely $p(e)$ of these circuits. 
\item{{\rm(ii)}}The graph $G$ does not contain the Petersen graph.

Thus it appears  useful to have a structural characterization
of graphs that do not contain the Petersen
graph, but that is undoubtedly a hard problem. In [\cite{RobSeyThoCubic}]
we managed to find such a characterization for cubic graphs 
under an additional connectivity assumption. We need a few definitions
before we can state the result. If $G$ is a graph
and $X\subseteq V(G)$, we denote by $\delta_G(X)$ or $\delta(X)$ the
set of edges of $G$ with one end in $X$ and the other in $V(G)-X$.
We say that a cubic graph is \dfn{theta-connected} if $G$ has girth
at least five,  and
$|\delta_G(X)|\ge6$ for all $X\subseteq V(G)$ such that
$|X|, |V(G)-X|\ge7$. We say that a graph $G$ is \dfn{apex} if
$G\backslash v$ is planar for some vertex $v$ of $G$ ($\backslash$ denotes
deletion). We say that a graph $G$ is \dfn{doublecross} if it has four
edges $e_1,e_2,e_3,e_4$ such that the graph  $G\backslash\{e_1,e_2,e_3,e_4\}$
can be drawn in the plane with the unbounded face bounded by a circuit
$C$, in which $u_1,u_2,v_1,v_2,u_3,u_4,v_3,v_4$ are pairwise distinct and occur
on $C$ in the order listed, 
where the edge $e_i$ has ends $u_i$ and $v_i$ for
$i=1,2,3,4$.
The graph \dfn{Starfish} is shown in Figure 1. Now we can state the
result of [\cite{RobSeyThoCubic}].

\goodbreak\midinsert
\vskip4.3in
\includegraphics{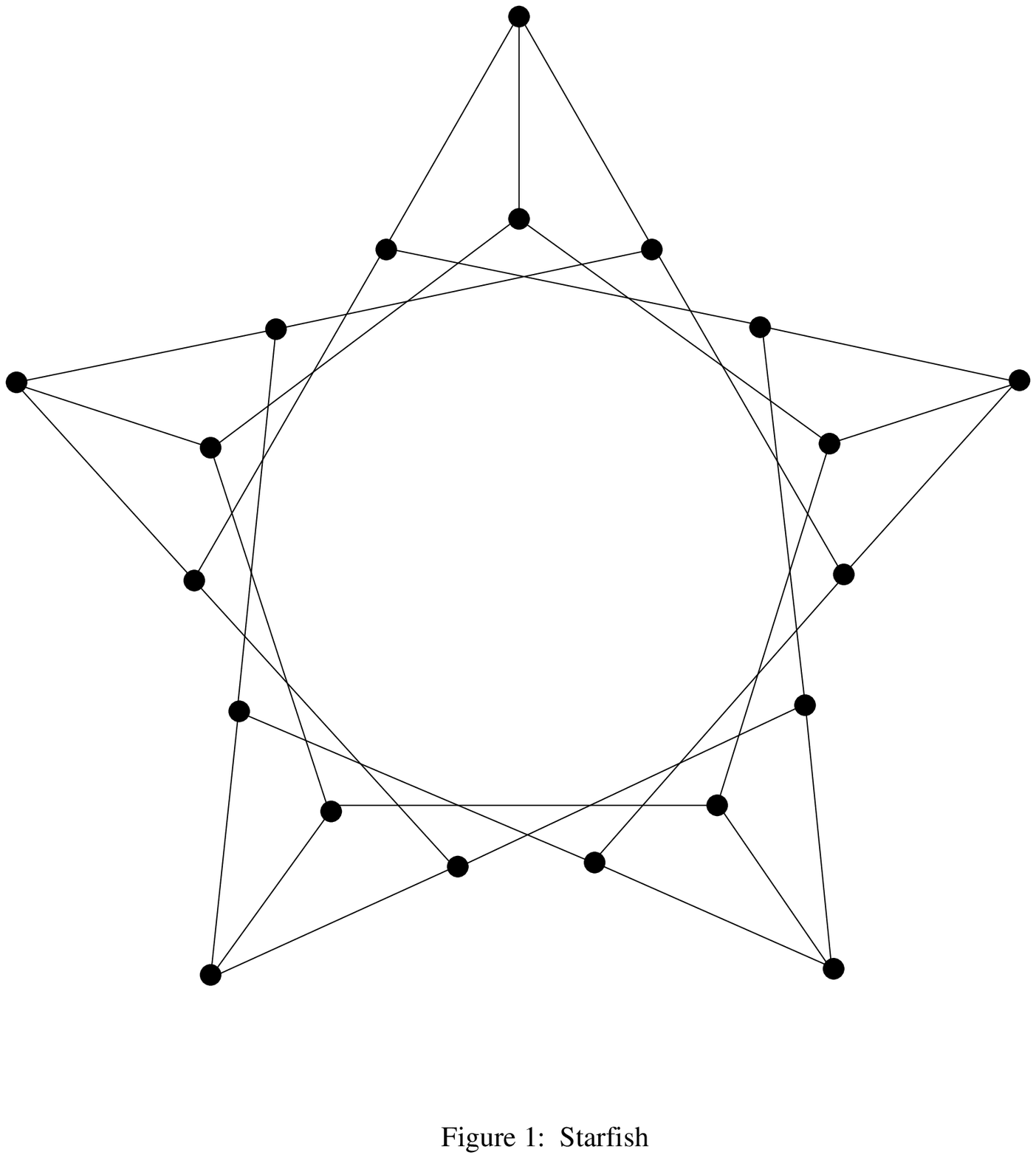}
\endinsert

\newthm{cubic}Let $G$ be a cubic theta-connected graph. Then $G$ does not contain
the Petersen graph if and only if either $G$ is apex, or $G$ is doublecross,
or $G$ is isomorphic to Starfish.

In the present paper we use \refthm{cubic} to prove the result stated
in the title and in the abstract, formally the following.

\newthm{main}Every cubic graph of girth at least six contains the Petersen
graph.

Theorem~\refthm{main} does not extend to graphs of minimum degree three.
For instance, let $H$ be (the $1$-skeleton of) the Dodecahedron.
The graph $H$ has an induced matching $M$ of size six.
Let $G$ be obtained from $H$ by subdividing every edge of $M$,
adding a new vertex $v$ and joining $v$ to all the vertices that resulted
from subdividing the edges of $M$.
Then $G$ has girth six, but it is apex, and hence does not contain the 
Petersen graph.

We prove \refthm{main} by induction, but in order for the inductive argument
to work we need to prove a stronger statement which we now introduce.
We say that two circuits of a graph \dfn{meet} if they have at least one
vertex in common. Thus if two circuits of a cubic graph meet, then they
have at least one edge in common.
We say that a circuit of a graph is \dfn{short} if it
has at most five edges. A short circuit of a graph $G$ which meets every
short circuit of $G$ is called a \dfn{breaker}.
We say that a graph $G$ is \dfn{interesting}
if it is cubic, it has at least ten vertices, and either it has girth at least
six or it has a breaker.
We shall see later that every interesting graph has at least fourteen
vertices. In fact, it can be shown that there is exactly one interesting
graph on fourteen vertices; this graph has girth six, and  is usually
called the \dfn{Heawood graph}.
 The result
we  prove is the following.

\newthm{main2}Every interesting cubic graph contains the Petersen graph.

\noindent Since the Petersen graph is not interesting, one might ask
if it is perhaps true that every interesting graph contains the Heawood
graph. Unfortunately, it is not. If a graph contains
another graph, and the first admits an embedding in the Klein Bottle
(in fact, any fixed surface), then so does the second. However, the
Heawood graph does
not admit an embedding in the Klein bottle, and yet there are cubic graphs
of girth six that do.

To prove \refthm{main2} we first show in Section~\refsec{ad}, using
\refthm{cubic}, that
\refthm{main2} holds for theta-connected interesting graphs, and
then prove \refthm{main2} for all interesting graphs
in Section~\refsec{int}. 

Andreas Huck (private communication) informed us that he can use \refthm{main}
to deduce the following result. A graph is \dfn{Eulerian} if every
vertex has even degree.

\newthm{huck}Let $G$ be a cubic $2$-edge-connected graph not containing the
Petersen graph. Then there exist five Eulerian subgraphs of $G$ such that
every edge of $G$ belongs to exactly two of these graphs.

\newsection{ad}APEX AND DOUBLECROSS GRAPHS

The objective of this section is to prove \refthm{theta} below,
our main theorem for theta-connected graphs.
We begin with the following.

\newthm{14verts}Every interesting cubic graph has at least
fourteen vertices.

\proof
Let $G$ be an interesting graph.
It is easy to see that every cubic graph of girth at least six  has
at least fourteen vertices. Thus we may assume that $G$ has a breaker $C$.
Let $H$ be the graph obtained from $G$ by deleting the edges of $C$.
Then $H$ has at most five vertices of degree one, and hence at least five
vertices of degree three, because $G$ has at least ten vertices. 
It follows that $H$ has a circuit. Let $C'$
be the shortest circuit of $H$.
The circuit $C'$ has length at least six, because it is disjoint from
$C$. Let $Z$ be the set of all vertices in $V(G)-V(C')$ that are adjacent
to a vertex of $C'$. Then $|Z|= |V(C')|$ by the choice of $C'$.
If $C'$ has length at least seven, then $|V(G)|\ge |V(C')|+|Z|\ge14$,
as desired, and so we may assume that $C'$ has length six. The
above argument shows that $G$ has at least twelve vertices, and so
we assume for a contradiction that $G$ has exactly twelve vertices. Thus
$V(G)=V(C')\cup Z$. Hence $V(C)\subseteq Z$, and the inclusion is
proper, because $C$ is short. 
The subgraph of $G$ induced by $Z-V(C)$ is $2$-regular, and hence has
a circuit.
But this circuit is short and disjoint from $C$,
a contradiction. Thus $G$ has at least fourteen vertices, as
desired.~\qed

%

A \dfn{pentagon} is a circuit of length five.

\newthm{2pent}Every two distinct pentagons in an interesting theta-connected
cubic graph have at most one edge in common.

\proof Let $G$ be an interesting theta-connected graph. Suppose
for a contradiction that $G$ has two distinct pentagons $C$
and $C'$ with more than one edge in common. Then
$|\delta_G(V(C)\cup V(C'))|= 5$ and $|V(C)\cup V(C')|= 7$, because
$G$ has girth at least five, and hence  $|V(G)-(V(C)\cup V(C'))|\le 6$
by the theta-connectivity of $G$, contrary to \refthm{14verts}.~\qed

If $G$ is a graph and $X\subseteq V(G)$, we denote by $G|X$ the graph
$G\backslash(V(G)-X)$.

\newthm{5pent}Every interesting theta-connected cubic graph has at 
most five pentagons.

\proof Let $G$ be an interesting theta-connected graph.
Since $G$ is theta-connected, every short circuit in $G$ is a pentagon.
Suppose for a contradiction that $G$ has at least six pentagons.  
Let $C_0$ be a breaker in $G$, and let
$C_1$, $C_2$, $C_3$, $C_4$ and $C_5$ be five other pentagons of $G$.
The sets $E(C_0)\cap E(C_i)$ ($i=1,2,\ldots,5$) are nonempty, and, by 
\refthm{2pent}, they are pairwise disjoint and each has cardinality one.
Thus $G$ has no other short circuit.
For $i=1,2,\ldots,5$ let $E(C_0)\cap E(C_i)=\{e_i\}$. We may assume
that $e_1,e_2,\ldots,e_5$ occur on $C_0$ in the order listed.
By \refthm{2pent} consecutive circuits in the sequence 
$C_1,C_2,\ldots,C_5, C_1$ have precisely one edge and its ends in common,
and non-consecutive circuits are vertex-disjoint, as otherwise $G$ has a short circuit distinct from $C_0,C_1,\ldots,C_5$. We conclude that 
$\left|\bigcup_{i=0}^5 V(C_i)\right|=15$.
Let $X=V(G)-\bigcup_{i=0}^5 V(C_i)$; then
$|\delta_G(X)| \le 5$, 
and hence $|X|\le 6$ by the theta-connectivity of $G$.
Moreover, $X$ has an odd number of elements, and hence is not empty.
Since $X$ has at most five vertices and is disjoint from $V(C_0)$, we
deduce that $G|X$ has no circuit, and
that $G|X$ is a path on at most three vertices.
Hence every vertex of $X$ is incident with an edge in $\delta(X)$,
and there exists a vertex $v\in X$ adjacent to every vertex
of $X-\{v\}$. 
We may assume $v$ has a neighbor $c_1 \in V(C_1)$, and some 
vertex $c_2 \in V(C_2)$ has a neighbor in $X$. Thus $c_1,c_2$ are joined by a two-edge 
path  with interior in $(C_1\cup C_2)\backslash V(C_0)$, and by a path of length at most three with 
interior in $X$, and their union is a short circuit disjoint from $C_0$, a contradiction.~\qed


\newthm{dblecross} Every  cubic doublecross graph of girth at least five 
has at least six pentagons.

\proof Let $G$ be a  doublecross graph of girth at least five, 
and let $e_1,u_1,v_1,\dots, e_4,u_4,v_4$ and $C$ be as in the
definition of doublecross. 
Let $P_1$ be the subpath of $C$ with ends
$u_1$ and $u_2$ not containing $v_1$,  let $P_2$ be the subpath of
$C$ with ends $u_2$ and $v_1$ not containing $v_2$, and let $P_3,P_4,
\dots, P_8$ be defined similarly.  Thus $C=P_1\cup P_2\cup \dots \cup P_8$,
and the paths $P_1,P_2,\dots, P_8$ appear on $C$ in the order listed.
Let $G':=G\backslash \{e_1,e_2,e_3,e_4\}$. We will regard $G'$ as
a plane graph with outer cycle $C$.
Let $f$ be the  number of bounded faces of $G'$,
and let $p$ be
the number of those that are bounded by a pentagon.
By Euler's formula
$|V(G)|+f+1=|E(G)|-4+2$, and since $G$ is cubic,
$2|E(G)|=3|V(G)|$. We deduce that $|E(G)|=3f+9$.
For $i=1,2,\dots, 8$ let $d_i=|E(P_i)|$. Since every edge of $G'\backslash E(C)$
is incident with two bounded faces, and every edge of $C$ is incident with one we obtain 
$2|E(G)|\ge 6f-p+8+\sum^8_{i=1} d_i$,
and hence $\sum^8_{i=1} d_i\le 10+p$.

Since $d_4,d_8\ge 1$, we get
$d_1+d_2+d_3+d_5+d_6+d_7\le 8+p$.
Let $q$ be the number of pentagons in the subgraph formed by $P_1,P_2,P_3$ 
and the edges $u_1v_1,u_2v_2$.
Then $d_1+d_2 + 1$ is at least five, and at least six unless the cycle 
$P_1\cup P_2+u_1v_1$ is a pentagon, and similarly for $d_2 + d_3 + 1$.
Furthermore,
$d_1 + d_3 +2$ is at least six unless $P_1\cup P_3 + u_1v_1 + u_2v_2$ is a pentagon. 
By adding,
$$(d_1+d_2 + 1) + (d_2 + d_3 + 1) + (d_1 + d_3 +2)
\ge 6 + 6 +6 -q;$$ 
that is, $d_1+d_2+d_3 \ge 7-q/2$.
Similarly if there are $r$ pentagons in the opposite crossing, then
$d_5+d_6+d_7 \ge 7 - r/2$.

So, adding,
$$8+p\ge d_1+d_2+d_3+d_5+d_6+d_7 \ge 14 - (q+r)/2.$$
So $p+(q+r)/2 \ge 6$, and hence $p+q+r\ge 6$ as required.~\qed

\newthm{apex}Every cubic apex graph of girth at least five has at least six pentagons.

\proof Let $G$ be a cubic apex graph of girth at least five,
and  let $v$ be a vertex of $G$ such that $G\backslash v$ is planar.
Let $f$ be the
number of faces in some planar embedding of $G\backslash v$,
and let $p$ be the number of them that are bounded by a pentagon. Then
$2|E(G)|=3|V(G)|$ and $|E(G)|=|E(G\backslash v)|+3$ because $G$ is cubic,  
$|V(G\backslash v)|+f=|E(G\backslash v)|+2$ by Euler's formula, and
$2|E(G\backslash v)|\ge 6f-p$,
since $G$ has girth at least five. We deduce that $p\ge 6$, 
as desired.~\qed 

In view of~\refthm{dblecross} and~\refthm{apex} it is natural
to ask whether every cubic graph of girth at least five not 
containing the Petersen graph has at least six pentagons. 
That is not true, because Starfish is a counterexample.

\newthm{theta}Every interesting theta-connected graph contains the Petersen
graph.

\proof Let $G$ be an interesting theta-connected graph. By \refthm{5pent}
$G$ has at most five pentagons. Thus by \refthm{dblecross}
$G$ is not doublecross, by \refthm{apex} $G$ is not apex, and $G$ is not
isomorphic to Starfish, because Starfish is not interesting. Thus $G$
contains the Petersen graph by \refthm{cubic}.~\qed

\newsection{int}INTERESTING GRAPHS

In this section we prove \refthm{main2}, which we restate below 
as \refthm{main3}.
Let $G$ be an interesting graph. We say that $G$ is \dfn{minimal} if
$G$ contains no interesting graph on fewer vertices.

\newthm{no3cut} Every minimal interesting graph has girth at least four. 

\proof 
Let $G$ be a minimal interesting graph; then $G$ is clearly connected.
Suppose for a contradiction that 
$C$ is a circuit in $G$ of length at most three.
Let $C'$ be a breaker in $G$,
and let $e\in E(C)\cap E(C')$. 
Let $H$ be obtained from $G$ by deleting $e$, deleting any
resulting vertex of degree one, and then suppressing all resulting vertices
of degree two. Then $H$ has at least ten vertices by \refthm{14verts}.
Also, it follows that
either $H$ has girth at least six or  $H$ has a breaker
(if $C\ne C'$ then the latter can be seen by 
considering the circuit of $H$ that corresponds to the circuit of
$C\cup C'\backslash e$). Thus
$H$ is interesting, contrary to the minimality of $G$.~\qed

We say that $X$ is a \dfn{shore} in a graph $G$ if $X$ is a set of vertices
of $G$ such that $|\delta(X)|\le5$ and
both $G|X$ and $G\backslash X$ have at least two circuits. 
The following is easy to see.

\newthm{shores}A cubic graph of girth at least five is 
theta-connected if and only if it has no shore.

\newthm{push} Let $G$ be an interesting graph, let $X$ be a
 shore in $G$,  and let $C$ be a breaker in $G$.
Then $G$ has a shore $Y$ such that $|\delta(Y)|\le|\delta(X)|$ and
$V(C)\cap Y=\emptyset$.

\proof Let $Y$ be a  shore in $G$ chosen so that 
$|\delta(Y)|$ is minimum, and subject to that,  $|Y\cap V(C)|$ is minimum.
We claim that $Y$ is as desired.
 From the minimality of $|\delta (Y)|$ we deduce that 
$|\delta (Y)|\le |\delta (X)|$ and that $\delta (Y)$ is a matching, 
and from the minimality of 
$|Y\cap V(C)|$ we deduce (by considering $V(G)-Y$) that $|Y\cap V(C)|\le 2$.
Suppose for a contradiction that $Y\cap V(C)\ne\emptyset$;
then $Y\cap V(C)$ consists of two vertices, say $u$ and $v$, that are
adjacent in $C$. Let
$Y'=Y-V(C)$. We deduce that $|\delta(Y')|\le
|\delta(Y)|$, and so it follows from the choice of $Y$ that $G|Y'$ has at most
one circuit. On the other hand since $u$ and $v$ are adjacent and have
degree two in $G|Y$ we see that $G|Y'$ has a circuit, and since
$|\delta(Y')|\le 5$ this is
a short circuit disjoint from $V(C)$, a contradiction.~\qed

\newthm{no5cut} No minimal interesting graph has a shore.


\proof Suppose for contradiction that $G$ is a minimal interesting graph, and
that $X$ is a shore in $G$ with $|\delta(X)|$ minimum. Let $k=|\delta(X)|$;
then $k\le 5$.
If $G$ has a short circuit let $C$ be a breaker in $G$;
otherwise let $C$ be the null graph.
By \refthm{push} we may assume that $V(C)\cap X=\emptyset$.
By \refthm{no3cut} and
the minimality of $k$ we may choose a circuit
$C'$ of $G\backslash X$ with $|V(C')|\ge k$. 
By the minimality of $k$ there exist $k$ disjoint
paths between $V(C')$ and $Z$, where $Z$ is the set of all vertices of
$X$ that are incident with an edge in $\delta(X)$. 
Let the paths be $P_1,P_2,\dots, P_k$, and for $i=1,2,\dots, k$ let
the ends of $P_i$ be $u_i\in Z$ and $v_i\in V(C')$ numbered so that
$v_1,v_2,\dots, v_k$ occur on $C'$ in this order. Let $H$ be obtained
from $G| X$ by adding a circuit $C''$ with vertex-set $\{w_1,w_2,\dots, w_k\}$
in order
and one edge with ends $u_i$ and $w_i$ for $i=1,2,\dots, k$.  Then $C''$
is a breaker in $H$.
Since
$G| X$ has a circuit, and that circuit, being disjoint from $C$, has
length at least six, we deduce that $H$ has at least ten vertices, and
so is interesting.  Moreover, $G$ 
contains  $H$, and is not isomorphic to $H$, because $G\backslash X$ is not
a circuit, contradicting the minimality of $G$.~\qed

We are now ready to prove \refthm{main2}, which we restate.

\newthm{main3} Every interesting graph contains the Petersen graph.

\proof It suffices to show that every minimal interesting graph contains
the Petersen graph.  To this end let $G$ be a minimal interesting graph.
If $G$ has girth at least five, then $G$ is theta-connected
by \refthm{shores} and \refthm{no5cut}, 
and hence contains the Petersen graph by
\refthm{theta}.  Thus we may assume that $G$ has a circuit of length
less than five, say $C$. Since $G$ has girth at least four
by \refthm{no3cut}, we deduce that $C$ has length four. 

We claim that $C$ is the only short circuit in $G$. To prove this claim
suppose for a contradiction that $C'$ is a short circuit in $G$ other
than $C$. Since $G$ is interesting,
we may assume that the pair $C,C'$ is chosen in such a way that
$C$ or $C'$ is a breaker in $G$.
Then $|\delta(V(C)\cup V(C'))|\le5$.
Let $X=V(G)-(V(C)\cup V(C'))$.
By \refthm{no5cut} $X$ is not a shore, and so $G|X$ has at most one circuit,
because $|\delta_G(X)|\le 5$. Thus $|X|\le5$.
It follows that $G$ has at most twelve vertices, contrary to
\refthm{14verts}.
Thus $C$ is the only short
circuit in $G$, as claimed.

Let the vertices of $C$ be $u_1,u_2,u_3,u_4$ (in order),  for
$i=1,2,3,4$ let $e_i$ be the edge of $C$ with ends $u_i$ and
$u_{i+1}$ (where $u_5$ means $u_1)$, and let $f_i$ be the unique
edge of $E(G)-E(C)$ incident with $u_i$. Let
$H$ be the graph obtained from $G\backslash e_1$ by contracting the edges
$e_2$ and $e_4$.  Then $H$ is a cubic graph with girth at least
five, and hence $|V(H)|\ge 10$, as is easily seen. Moreover, every
pentagon in $H$ contains one end of $e_3$.

Let us assume first that $H$ is theta-connected. By the minimality of
$G$, the graph $H$ is not interesting; in particular, the edge $e_3$
belongs to no pentagon of $H$. Since no two pentagons in a cubic graph
of girth at least five share more than two edges, we deduce that 
the ends of $e_3$
belong to at most two pentagons each. Thus $H$ is not doublecross
by \refthm{dblecross}, it is not apex by \refthm{apex}, and it is not isomorphic
to Starfish, because Starfish has three pairwise vertex-disjoint pentagons.
Thus $H$ contains the Petersen graph by
\refthm{cubic}, and hence so does $G$, as desired.  

We may therefore assume that $H$ is not theta-connected. By \refthm{shores}
$H$ has a shore. By \refthm{no5cut} applied to $G$
there exists a set $X_1\subseteq V(G)$  such that $|X_1|,|V(G)-X_1| \ge7$,
$|\delta_G (X_1)|=6$,  
$u_1,u_4\in X_1$, $u_2,u_3\not\in X_1$,
and that $\delta_G(X_1)$ is a matching. 
Thus $|X_1|, |V(G)-X_1|\ge 8$. 
By arguing similarly for the graph $G\backslash e_2$
we deduce that either $G$ contains the Petersen graph, or
there exists a set $X_2\subseteq V(G)$ such that 
$|X_2|, |V(G)-X_2|\ge 8$,
$|\delta_G (X_2)|=6$, $u_1,u_2\in X_2$, and $u_3,u_4\in V(G)-X_2$. We may
assume the latter.
Since $|\delta_G (X_1\cap X_2)|+|\delta_G (X_1\cup X_2)|\le |\delta_G (X_1)|+
|\delta_G(X_2)|=12$, we deduce that  $\delta_G (X_1\cap X_2)$ or
$\delta_G (X_1\cup X_2)$ has at most six elements. From the symmetry we may
assume that it is the former.  Since $e_1,e_4\in \delta_G(X_1\cap X_2)$,
it follows that $|\delta_G (Y)|\le5$, where
$Y=X_1\cap X_2-\{u_1\}$. By \refthm{no5cut} $Y$ is not a shore, and hence
the graph $G|Y$ has at most one circuit. However, if $G|Y$ has a circuit,
then that circuit does not meet $C$, and yet it has length at most
five (because $|\delta(Y)|\le5$), which is impossible. Thus $G|Y$ has
no circuit, and hence $|Y|\le3$.
Similarly, either $|X_1-X_2-\{u_4\}|\le 3$ or $|X_2-X_1-\{u_2\}|\le 3$,
and from the symmetry we may assume the former. Thus $|X_1|\le 8$.
Since $|X_1|\ge 8$ as we have seen earlier, the above inequalities
are satisfied with equality.  In particular, $|\delta_G (X_1\cup X_2)|=6$,
and hence
$|\delta_G (X_1\cup X_2\cup \{u_3\})|\le 5$, and likewise
$|\delta_G (X_2- X_1- \{u_2\})|\le 5$.
As above we deduce that $|V(G)-X_1|=8$. Thus
$G| X_1$ and $G\backslash X_1$ both have eight vertices, six vertices of
degree two, two vertices of degree three, and girth at least six. It follows
that $G| X_1$ and $G\backslash X_1$ are both isomorphic to the graph that
is  the union of three paths on four vertices each, with the same ends
and otherwise vertex-disjoint.

We now show that $H$ is isomorphic to the graph shown in Figure~2.
Let $G|X_1$ consist of three paths $ab_ic_id$ for $i = 1, 2,3$, and let
$G\backslash X_1$ have three paths $pq_ir_is$ similarly. 
So there is a six-edge
matching $M$ in $G$ between $\{b_1,b_2,b_3,c_1,c_2,c_3\}$ and 
$\{q_1,q_2,q_3,r_1,r_2,r_3\}$. 
Since $C$ exists we can assume that $b_3$ is matched by $M$ to $q_3$ and 
$c_3$ to $r_3$. 
Thus $V(C)=\{b_3,c_3,r_3,q_3\}$.

Now $X_2$ exists and contains $a,p$ and not $d,s$; and each of the six paths 
of the previous paragraph includes an edge of $\delta(X_2)$.
Thus no edge of $M$ belongs to $\delta(X_2)$. If $X_2$ 
contains both $b_1$ and $c_1$, then
$b_1,c_1$ are matched by $M$ to vertices with distance at most two in $X_2$, 
and hence $G$ has a circuit of length at most five disjoint from $C$, a contradiction. 
Thus $X_2$ contains at most one of $b_1,c_1$, and similarly for the pairs
$(b_2,c_2)$,  $(q_1,r_1)$ and $(q_2,r_2)$. We may therefore assume that 
$X_2 = \{a,b_1,b_2,b_3,p,q_1,q_2,q_3\}$. 
So $b_1,b_2$ are matched to $q_1,q_2$ and $c_1,c_2$ to $r_1,r_2$ in some order. 
Because $G$ is interesting, we may assume the pairs are
$(b_1,r_1)$, $(b_2,r_2)$, $(c_1,r_2)$, $(c_2,r_1)$.
Thus $H$ is isomorphic to the graph shown in Figure~2.
That graph, however, contains the Petersen graph, as desired.\qed

\goodbreak\midinsert
\vskip3in
\includegraphics{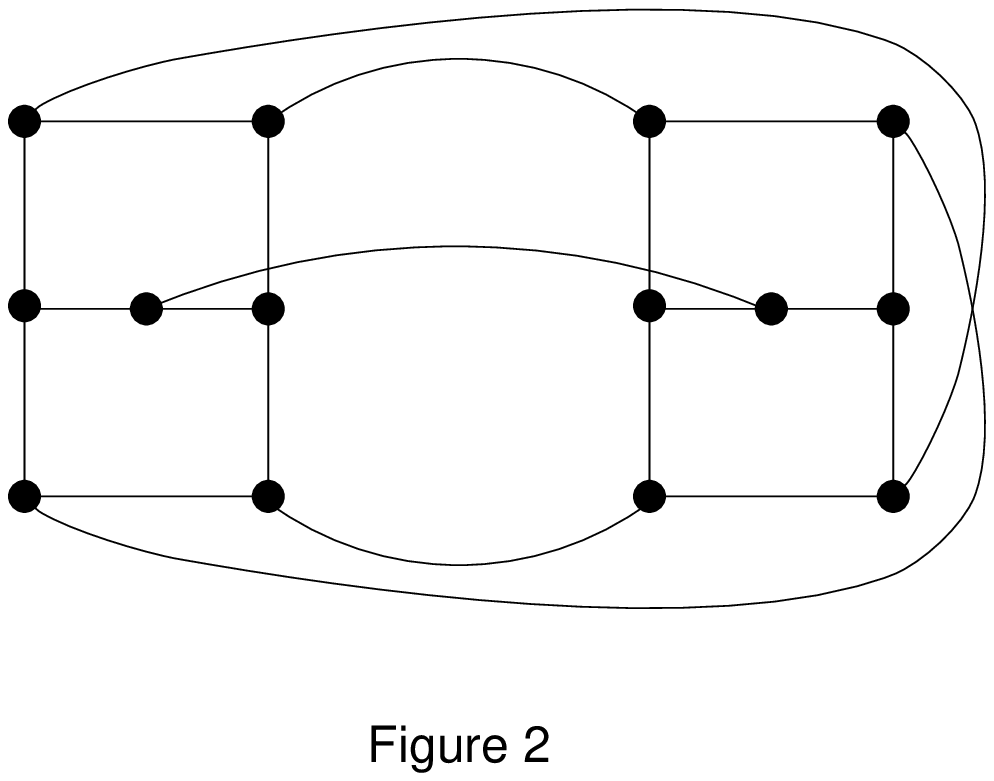}
\endinsert

\beginsection References


\def\TAMS{{\it Trans.\ Amer.\ Math.\ Soc.}}

\bibitem{AlsGodZha}B.~Alspach, L.~Goddyn and C.-Q.~Zhang,
Graphs with the circuit cover property,
\TAMS\ {\bf344} (1994), 131--154.

\bibitem{RobSeyThoCubic} N.~Robertson, P.~D.~Seymour and R.~Thomas,
Excluded minors in cubic graphs, {\tt arXiv:1403.2118}.

\vfill

\noindent
This material is based upon work supported by the National Science Foundation.
Any opinions, findings, and conclusions or recommendations expressed in
this material are those of the authors and do not necessarily reflect
the views of the National Science Foundation.

\end

We say that $X$ is a \dfn{shore} if $X$ is a non-empty proper subset of $V(G)$,
and there do not exist sets $X_1,X_2\subseteq V(G)$ such that
$\delta (X)=\delta (X_1)\cup \delta (X_2)$, 
$\delta (X_1)\cap \delta (X_2)=\emptyset$
and $\delta (X_1)\ne\emptyset\ne \delta (X_2)$. This is equivalent to 
$k(G| X)+k(G\backslash X)=k(G)+1$, where $k(H)$ denotes the number of
components of $H$ and $G|X$ denotes the graph $G\backslash(V(G)-X)$. 
We say that $X$ is a \dfn{cyclic shore} if $X$ is
a shore, and both $G|X$ and $G\backslash X$ have at least one circuit.
We say that $X$ is a \dfn{nontrivial cyclic shore} if $X$ is a 
cyclic shore and neither $G| X$ nor $G\backslash X$ is a circuit.
Thus a cubic graph of girth at least five is theta-connected if and only if
it has no nontrivial cyclic shore $X$ with $|\delta(X)|\le5$.